\documentclass{gtart}

\input gtoutput
\volumenumber{3}\papernumber{13}\volumeyear{1999}
\pagenumbers{303}{330}\published{13 September 1999}
\proposed{Steve Ferry}
\seconded{Walter Neumann, David Gabai}
\received{19 December 1998}
\accepted{27 August 1999}

\usepackage{amssymb}

\newtheorem{Lemma}{Lemma} 
\newtheorem{Prop}{Proposition}
 
\newtheorem*{Sublemma}{Sublemma}
\newtheorem*{TA}{Theorem A}
\newtheorem*{TB}{Theorem B}
\newtheorem*{TC}{Theorem C}

\theoremstyle{remark}
\newtheorem*{Remark}{Remark}
\newtheorem{Example}{Example}
\newtheorem*{QA}{Question A}
\newtheorem*{QB}{Question B}

\begin{document}

\title{Piecewise Euclidean Structures and Eberlein's\\Rigidity Theorem
in the Singular Case}
\shorttitle{Piecewise Euclidean Structures and Eberlein Rigidity}

\author{Michael W Davis\\Boris Okun\\Fangyang Zheng}
\shortauthors{Davis, Okun and Zheng}

\asciiaddress{Department of Mathematics, The Ohio
State University\\ Columbus, OH 43201, USA\\
Department of Mathematics, Vanderbilt University\\ Nashville, TN
                                  37400, USA\\
Department of Mathematics, The Ohio
State University\\ Columbus, OH 43201, USA
}

\asciiemail{mdavis@math.ohio-state.edu, okun@math.vanderbilt.edu, zheng@math.ohio-state.edu}

\address{Department of Mathematics, The Ohio
State University\\ Columbus, OH 43201, USA\\\smallskip
Department of Mathematics, Vanderbilt University\\ Nashville, TN
                                  37400, USA\\\smallskip
Department of Mathematics, The Ohio
State University\\ Columbus, OH 43201, USA\\\smallskip
{\rm Email:}\qua\tt mdavis@math.ohio-state.edu, okun@math.vanderbilt.edu\\
zheng@math.ohio-state.edu}

\begin{abstract}
In this article, we generalize Eberlein's Rigidity Theorem to the
singular case, namely, one of the spaces is only assumed to be a
CAT(0) topological manifold. As a corollary, we get that any compact 
irreducible but
locally reducible locally symmetric space of noncompact type does not
admit a nonpositively curved (in the Aleksandrov sense) piecewise
Euclidean structure. Any hyperbolic manifold, on the other hand, does
admit such a structure.
\end{abstract}



\maketitlepage

\section{Introduction}

Suppose that $V$ and $V^{\ast }$ are compact locally symmetric
manifolds with fundamental groups $\Gamma$ and $\Gamma^{\ast }$ and
with universal covers $X$ and $X^{\ast }$, respectively. Suppose
further that the symmetric spaces $X$ and $X^{\ast }$ have no compact
or Euclidean factors and that no finite cover of $V^{\ast }$ splits
off a hyperbolic 2--manifold as a direct factor. Let $\theta\co \Gamma
\to \Gamma^{\ast }$ be an isomorphism.  The Mostow Rigidity Theorem
asserts that the isomorphism $\theta$ is induced by an isometry $V \to
V^{\ast }$, possibly after renormalizing the metric on $V^{\ast }$.
(``Renormalizing the metric'' means that we are allowed to rescale the
metric on each factor of $X^{\ast }$ by multiplying by a positive
constant.)

In the early eighties, Eberlein, Gromov and Heintze independently
proved generalizations of Mostow's Theorem in which the hypotheses on
$V$ are weakened.  They showed that the conclusion of Mostow's Theorem
remains valid without the assumption that $V$ is locally symmetric.
It is only necessary to assume that it is a closed Riemannian manifold
of nonpositive sectional curvature, provided that the rank of $X^{\ast
}$ is at least two.  (If the rank of $X^{\ast }$ is one, there is no
such generalization. For, if $X^{\ast }$ is of rank one, then the
sectional curvature of $V^{\ast }$ is strictly negative and hence, the
same will be true for any sufficiently small deformation of its
metric; moreover, the new metric will not, in general, be isometric to
the original one.)

Eberlein proved the above generalization of the Rigidity Theorem in
the case where $X^{\ast }$ is reducible and the uniform lattice
$\Gamma^{\ast }$ is irreducible.  We shall sometimes refer to this as
``Eberlein's case''.  Gromov proved the result in full generality;
however, the argument has only appeared (in \cite{BGS}) in the case
where $X^{\ast }$ is irreducible.

In this paper we will prove a further generalization of Eberlein's
case of the Rigidity Theorem by allowing a further weakening of the
hypotheses on $V$: the metric need not be smooth. To explain this, we
first need to understand what ``nonpositive curvature'' means in this
context.

About fifty years ago Aleksandrov showed that the notions of upper or
lower curvature bounds make sense for a more general class of metric
spaces than Riemannian manifolds, namely, they make sense for
``geodesic spaces''.  A {\em geodesic segment\/} in a metric space $Y$
is the image of an isometric embedding of an interval into $Y$; the
space $Y$ is a {\em geodesic space\/} if any two points can be connected
by a geodesic segment.  A {\em triangle\/} in $Y$ is a configuration of
three points (the vertices) and three geodesic segments (the edges)
connecting them.  Following Aleksandrov, one defines $Y$ to be {\em
nonpositively curved\/} if the distance between any two points in any
small triangle in $Y$ is no greater than the distance between the
corresponding points of a comparison triangle in the Euclidean
plane. (The precise definition will be given in Section 1.1, below.)

The weakened hypotheses on $V$ will then be that it is a complete geodesic
space, that it is nonpositively curved in
the above sense and that it is a closed topological manifold.

Our interest in generalizing the Rigidity Theorem to the singular case
arose from an attempt to answer the following basic question
concerning nonpositively curved, polyhedral metrics on manifolds.
Given a Riemannian manifold $M$ of nonpositive sectional curvature, is
it possible to find a piecewise Euclidean, polyhedral metric on $M$
which is nonpositively curved in Aleksandrov's sense?  An application
of our rigidity result is that this is, in fact, impossible when
$M=V^{\ast }$. Indeed, if $V$ were a nonpositively curved polyhedron
homeomorphic to $V^{\ast }$, then by the Rigidity Theorem it would be
isometric to $V^{\ast }$, which is clearly impossible (since $V$
contains open subsets which are flat and $V^{\ast }$ does not). A
further discussion of the problem of approximating nonpositively
curved Riemannian metrics by nonpositively curved polyhedral ones will
be given in Section 1.

Next we shall make a few remarks concerning the proofs of Mostow and
Eberlein.  Since $V$ and $V^{\ast }$ are Eilenberg--MacLane spaces, the
isomorphism $\theta \co \Gamma \to\Gamma^{\ast }$ gives rise to a
homotopy equivalence $\overline{f}\co V\to V^{\ast }$. This lifts to a
$\theta$--equivariant map $f\co X\to X^{\ast }$. By analyzing the
maximal flats in $X$ and $X^{\ast }$, Mostow first shows that $f$
extends to a homeomorphism of the maximal (or Furstenberg)
boundaries. (The maximal boundary of $X^{\ast }$ is a homogeneous
space of the form $G/P$, where $G$ is the identity component of the
isometry group $I(X^{\ast })$ and $P$ is a minimal parabolic
subgroup.)  Let $X^{\ast }=X_1^{\ast } \times \cdots \times X_k^{\ast
}$ be the de Rham decomposition into irreducible factors so that
$I(X^{\ast })=I(X_1^{\ast }) \times \cdots \times I(X_k^{\ast })$.
The proof of Mostow's Theorem is the easiest to finish in Eberlein's
case (where $X^{\ast }$ is reducible and $\Gamma^{\ast }$ is
irreducible).  The key ingredient is the Borel Density Theorem
(cf \cite{R}).  It implies that the image of $\Gamma^{\ast }$ is
dense in the identity component of each $I(X_i^{\ast })$.  It follows
that the map of maximal boundaries is actually induced by an
isomorphism of Lie groups $I_0(X) \to I_0(X^{\ast })$ (where $I_0$
denotes the identity component of the isometry group).  From this, the
result follows (cf Corollary 18.2, page 133, in \cite{M}). (A
historical remark: one of the earliest rigidity results was proved by
Selberg in 1957.  He considered the case where $X^{\ast }$ is the
product of $k>1$ copies of the hyperbolic plane and where
$\Gamma^{\ast }$ is irreducible. He used the fact that $\Gamma^{\ast
}$ was dense in each copy of $I_0({\mathbb H}^2)=PSL(2,{\mathbb R})$
to show that any deformation of $\Gamma^{\ast }$ in $I(X^{\ast })$ had
to be trivial.  The details of this argument can be found in
\cite{S}.)

There are two steps in Eberlein's generalization of Mostow's argument
(to the case where $X$ is not required to be symmetric). First one
shows that $X$ also splits isometrically as $X=X_1 \times \cdots
\times X_k$.  (This is accomplished in \cite{E1}.)  In the second step
(which was actually accomplished in the earlier papers \cite{CE} and
\cite{E2}) one shows that each factor $I_0(X_i)$ is actually a simple
Lie group and that each $X_i$ is a symmetric space. Hence, one is
reduced back to Mostow's Theorem.  The arguments in the second step
are mainly Lie group theoretic. The key ingredient is again a density
property (which is closely connected to the Borel Density Theorem),
namely, that $\Gamma$ satisfies the ``duality condition'' of
\cite{CE}.  One shows that the the image of $\Gamma$ in $I(X_i)$ is
not discrete and hence, that the identity component $G_i$ of the
closure of its image in $I(X_i)$ is a nontrivial Lie group.  One then
argues that each $G_i$ must be semisimple since otherwise $\Gamma$
(and hence $\Gamma^{\ast }$) would contain a nontrivial normal abelian
subgroup.

Our argument, in the case where $X$ is no longer assumed to be
Riemannian, follows the same general outline.  We show, in Section
2.1, that Eberlein's arguments can be modified to prove that $X$
splits isometrically, as before.  As for the second step, we first
must face the fact that we don't know, {\em a priori,} that the
isometry groups $I(X_i)$ are Lie groups.  (The Hilbert--Smith
Conjecture is still open --- if a topological group acts effectively on
a topological manifold it is not known if this forces it to be a Lie
group.)  In fact, {\em a priori,} the spaces $X_i$, which are manifold
factors, need not be manifolds.  We get around this by showing, in
Section 2.3, that $I(X_i)$ acts transitively on $X_i$. (This comes
from the fact that $\Gamma$ is dense in $I_0(X_i)$.)  It then follows
from a theorem of Montgomery and Zippin that $I(X_i)$ is a Lie group.
Then, following the line of argument in \cite{CE} and \cite{E2} one
can conclude that each $I_0(X_i)$ is a simple Lie group.  At this
stage we are still not done.  Even though $X$ is now a homogeneous
space of the form $G/K$, the possibility remains that its
$G$--invariant metric is not Riemannian.  However, we show in Lemma 18
that the nonpositivity of the curvature excludes this possibility.

After writing a preliminary version of this paper we learned of
B Leeb's recent preprint \cite{Le}.  A corollary of his theorem is
that the Rigidity Theorem also holds in the singular case whenever
each factor of $X^{\ast }$ has rank $\geq 2$.  Leeb's approach
generalizes Gromov's line of argument in \cite{BGS}.  Thus, the
singular version of the Rigidity Theorem holds in complete generality.
We should also mention that Bruce Kleiner can prove our result by a
somewhat different argument.

Davis and Zheng were partially supported by
NSF grants, Zheng was also partially supported by an 
Alfred P Sloan Fellowship.  This project is also sponsored by the 
National Security Agency under grant \# MDA 904-98-1-0036.

%%\tableofcontents

\section{ Piecewise Euclidean CAT(0) structures }
\subsection{Hadamard spaces}
A metric space is a {\em length space\/} if the distance between any two
points is the infimum of the lengths of all paths between them.  The
space is a {\em geodesic space\/} if the infimum is always realized by a
path of minimum length.

Given any triangle $\Delta $ in a geodesic space $X$ there is a
triangle $\Delta^{\ast }$ in the Euclidean plane with the same edge
lengths and a well defined isometry $ \Delta \to \Delta^{\ast }$ which
is denoted by $x\to x^{\ast }$. The CAT(0)--inequality of \cite{G}
asserts that for any two points $x,y\in \Delta$, $d(x,y)\leq |x^{\ast
}-y^{\ast }|$. (Thus, any triangle in $X$ is ``thinner'' than the
corresponding triangle in ${\mathbb R}^2$).

A complete geodesic space $X$ is a {\em Hadamard space\/} (or a {\em CAT(0)
space\/}) if each triangle satisfies the CAT(0)--inequality.

A geodesic space is {\em nonpositively curved\/}, if the
CAT(0)--inequality holds locally. It can be proved (cf \cite{B} or
\cite{BH}) that if a geodesic space is complete and nonpositively
curved, then its universal cover is a Hadamard space.

 \begin{Example}[(Aleksandrov and Toponogov)]
A complete, simply connected \linebreak Riemannian manifold is a Hadamard space
if and only if its sectional curvature is $\leq 0$.
\end{Example}

A key fact is that, for any Hadamard space $X$, the metric $d\co X
\times X \to {\mathbb R}$ is a convex function, cf \cite{G}. (This
means that its restriction to any geodesic segment in $X \times X$ is
a convex function.)  It follows that any Hadamard space is
contractible.  (There is an unique geodesic segment from a given base
point to any other point; one can then contract the space to the base
point by shrinking the geodesic segments.)

A geodesic space has {\em extendible geodesics\/} if every geodesic
 segment can be extended to a larger interval.  A geodesic space
is {\em geodesically complete\/} if every geodesic segment can
be extended to a geodesic line. An easy argument shows that if a complete 
geodesic space  has extendible geodesics then
it is geodesically complete.

Suppose that the geodesic space $X$ is nonpositively curved.
If $X$ is a topological manifold (or even if its local homology $H_*(X, X-x)$
is nontrivial at each point $x \in X$), then it has extendible geodesics.
Indeed, suppose $c\co [0,d] \to X$ is a geodesic segment which does not
extend past $c(0)$.  Choose a small ball $B_{\epsilon}$ centered at $c(0)$
with $\epsilon < d$ and small enough so that $B_{\epsilon}$ is convex.
Since every point in $B_{\epsilon} - c(0)$ can be joined by a geodesic
segment to $c(\epsilon )$, $B_{\epsilon} - c(0)$ contracts to $c(\epsilon)$,
contradicting the assumption that $H_*(B_{\epsilon}, B_{\epsilon} - c(0))$ is
nonzero.

Given two subsets $A$ and $B$ of $X$, the {\em Hausdorff distance\/} 
between them,
denoted by $Hd(A,B)$, is the smallest number $\epsilon$ such that each
subset is contained within an $\epsilon$--neighborhood of the other.
The subsets $A$ and $B$ are {\em parallel\/} if their Hausdorff distance
is finite.  An {\em $r$--flat\/} in $X$ is a subset which is isometric to
$r$--dimensional Euclidean space.

A subset $A$ of a geodesic space is {\em locally convex\/} if the
geodesic segment between any two sufficiently close points in $A$
actually lies in $A$.  The subset $A$ is {\em totally geodesic\/} if any
geodesic line which contains two distinct points of $A$ is contained
in $A$.

Given a closed convex subset $A$ in a Hadamard space $X$, there is a
continuous {\em projection map} $p\co X \to A$ which sends a point $x
\in X$ to the nearest point in $A$.

If the Hadamard space $X$ is locally compact, then it can be
compactified by adding an ``ideal boundary'' $X(\infty )$. One way to
define this is as the set of all geodesic rays emanating from a given
basepoint.  If $X$ is Riemannian, then $X(\infty )$ is homeomorphic to
a sphere, but in the general case it can be much more complicated.

Let $\gamma$ be an isometry of $X$.  Its {\em translation distance\/}
$d(\gamma)$ is the infimum over all points $x \in X$ of $d(x,\gamma
x)$.  The isometry $\gamma$ is {\em semisimple\/} if there is at least
one point $x$ such that $d(\gamma)=d(x,\gamma x)$; $\gamma$ is a {\em
Clifford translation\/} if $d(\gamma)=d(x,\gamma x)$ for all $x \in X$.
If $\gamma$ is semisimple, then its {\em minimum set} $MIN(\gamma)$ is
the set of $x$ in $X$ such that $d(\gamma )=d(x,\gamma x)$.  When $X$
is Hadamard, it follows from the convexity of the metric that
$MIN(\gamma)$ is a convex subset.

A length space is {\em singular\/} if it is not isometric to a
Riemannian manifold. An important class of singular metrics is the
class of polyhedral metrics, defined below.

Suppose that $K$ is a cell complex formed by gluing together convex
polytopes in some Euclidean space via isometries between certain
faces. Then $K$ is said to have a {\em piecewise Euclidean
structure.} For example, if the cell complex $K$ is a piecewise
linearly embedded subcomplex of Euclidean space (ie, if it is a cell
complex in the classical sense), then it has such a structure. A
piecewise Euclidean structure on $K$ defines a length metric: the
distance between two points is the infimum of the lengths of all
piecewise linear paths between them. Such metrics are called {\em
piecewise Euclidean\/} (or {\em polyhedral\/}).

For the further details concerning the definitions and basic
properties of CAT(0) spaces, we refer the reader to the excellent
books of Bridson--Haefliger \cite{BH} and Ballmann \cite{B}.

\subsection{Convex hypersurfaces}
It follows from the Gauss equation that a smooth, convex hypersurface
in the Euclidean space has sectional curvature $\geq 0$. The singular
analogue of this is the following. The boundary $K$ of a convex
polytope in Euclidean space has a piecewise Euclidean metric and this
metric is nonnegatively curved in the sense that the reverse
inequality to the CAT(0)--inequality is satisfied. In fact, it was
precisely for the purpose of studying polyhedral metrics on convex
polyhedral surfaces that Aleksandrov initiated the study of length
spaces with curvature bounded from below.

The Gauss equation remains valid for the hypersurfaces in manifolds
equipped with an indefinite metric, the only change that is necessary
is that one side of the equation must be multiplied by $-1$ when the
inner product is negative definite in the normal direction. As we
shall see, this sign change allows the possibility that a convex
hypersurface in Minkowski space will be nonpositively curved.

Recall that {\em Minkowski space} ${\mathbb R}^{n,1}$ is an
$(n+1)$--dimensional real vector space with coordinates $(x_0, x_1,
\ldots , x_n)$ equipped with a {\em Lorentzian\/} inner product of
signature $(n,1)$ (eg, $x\cdot x = -x_0^2 + x_1^2 + \cdots +
x_n^2$). A hyperplane in $\ {\mathbb R}^{n,1}$ is {\em spacelike\/} if
its normal vector $n$ is {\em timelike}, ie, if $n\cdot n<0$. A
smooth hypersurface in ${\mathbb R}^{n,1}$ is {\em spacelike\/} if its
tangent space at each point is spacelike. Since the restriction of the
inner product to any spacelike hyperplane is positive definite, any
spacelike hypersurface inherits a Riemannian structure. By the Gauss
equation any smooth, spacelike convex hypersurface in ${\mathbb
R}^{n,1}$ is nonpositively curved.

By definition a {\em convex polyhedral hypersurface} $\Sigma $ in $\
{\mathbb R}^{n,1}$ is the boundary of an $(n+1)$--dimensional convex
polyhedral set $C$ in ${\mathbb R}^{n,1}$ (ie, $C$ is defined by a
discrete set of affine inequalities). The convex hypersurface $\Sigma
$ is {\em spacelike\/} if each supporting hyperplane of $C$ is
spacelike. As in the Euclidean case, $\Sigma $ has an induced
piecewise Euclidean metric.

 \begin{Example} [\cite{CDM}] 
Any spacelike convex polyhedral hypersurface in ${\mathbb R}^{n,1}$ is
a Hadamard space.  (Recently, Moussong has extended this result to
arbitrary spacelike convex hypersurfaces in ${\mathbb R}^{n,1}$.)
\end{Example}

As a simple corollary we get the following example.

 \begin{Example} [\cite{CDM}]  Every complete hyperbolic manifold
$M^n$ admits a nonpositively curved, piecewise Euclidean metric. The
construction goes as follows. Identify $M^n$ with ${\mathbb
H}^n/\Gamma $ where $\Gamma $ is a discrete torsion-free subgroup of
$O_0(n,1)$, the isometry group of ${\mathbb H}^n$. Choose a net $S$ in
$M^n$ (ie, $S$ is a discrete subset of $M^n$ such that each point of
$M$ lies within a bounded distance of $S$) and let $\widetilde{S}$
denote the inverse image of $S$ in ${\mathbb H}^n$. Let $\Sigma $
denote the boundary of the convex hull of $\widetilde{S}$ in ${\mathbb
R}^{n,1}$. Then $\Sigma $ is a $\Gamma $--stable, spacelike, convex
polyhedral hypersurface in ${\mathbb R}^{n,1}$. Hence, as in Example
2, it is CAT(0). Thus, $\Sigma /\Gamma $ gives the desired piecewise
Euclidean metric on $M^n$. (The details of this construction can be
found in \cite{CDM}.)  Similarly, the product of two hyperbolic
manifolds can be given a piecewise Euclidean structure.  However, if,
for example, $\Gamma$ is an irreducible uniform lattice in $I({\mathbb
H}^2 \times {\mathbb H}^2)$, then, as we shall see in the next
section, the resulting 4--manifold is rigid.
\end{Example}

\begin{Remark} Gromov has suggested that the construction in Example 3
can be modified to show that any Riemannian manifold of sufficiently
pinched negative sectional curvature can be given a nonpositively
curved, piecewise Euclidean metric.
\end{Remark}

\subsection{Questions and further examples}
At this point we have two classes of examples of nonpositively curved
geodesic metrics on manifolds: 1) Riemannian metrics and 2) piecewise
Euclidean metrics on polyhedral manifolds. What is the intersection of
these two classes? When can one type of metric be deformed (through
nonpositively curved metrics) into the other type? Both classes
contain the Riemannian manifolds of constant curvature $0$ (ie, flat
manifolds). Example 3 shows that a complete metric of constant
curvature $-1$ can be deformed to a polyhedral metric of nonpositive
curvature. Before continuing the discussion let us break our problem
into two questions.

 \begin{QA} Suppose $M^n$ is a nonpositively curved Riemannian
manifold. Does $M^n$ admit a nonpositively curved piecewise
Euclidean metric?
\end{QA}

\begin{QB} Suppose $M^n$ is a nonpositively curved piecewise
Euclidean polyhedral manifold. Does $M^n$ admit a Riemannian metric of
nonpositive sectional curvature?
\end{QB}

In dimension 2 both questions have affirmative answers. In dimension 3
the answers are not known (although it seems likely that the answer is
``yes" in both cases). As we shall see in Example~6, below, in
dimension $\geq 5$, Question B has a negative answer, while in
dimension 4 the answer is unknown.  A corollary to the main result of
this paper is that, in general, Question A also has a negative answer
in dimensions $\geq 4$.  As we have seen in Example~3, Question~A
has a positive answer when $M^n$ is a real hyperbolic manifold.  The
authors speculate that the answer will be negative for the other
locally symmetric manifolds, eg for complex and quaternionic
hyperbolic manifolds.

 \begin{Example}[(\cite{G}, \cite{CD1})]  Suppose that $V$ is a
manifold with a (possibly singular) length metric, that $W$ is a
codimension-two submanifold and that $V_k$ is a $k$--fold cyclic
branched cover of $V$ along $W$.  (One way to insure that such
branched cover exist, for any $k$, is to require $V$ and $W$ to be
closed and oriented and the image $[W]$ of the orientation class to be
zero in $H_{n-2}(V)$.)
\end{Example}

\begin{TA}   Let $V$, $W$ and $V_k$ be as the above.
For the induced metric on $V_k$ to be nonpositively curved
it is necessary and sufficient that the following two conditions hold:

{\rm(a)}\qua $V$ is nonpositively curved, and

{\rm(b)}\qua $W$ is a locally convex subset of $V$.
\end{TA}

\begin{proof} The necessity of the conditions (a) and (b) is obvious.
As for the sufficiency, suppose $\widetilde{\Delta }$ is a small
geodesic triangle in $V_k$ and that $\Delta $ denotes its image in
$V$.  If $\widetilde{\Delta} $ intersects the branch set $W$ in a
vertex or an edge or if it is contained in $W$, then condition (a) and
(b) insure that the distance between two points of $\widetilde{\Delta
}$ is the same as that of the corresponding points of $\Delta $. If
the intersection of $\widetilde{\Delta}$ and $W$ is one or more points
in the interiors of the edges, then we can subdivide $\widetilde
{\Delta}$ and apply a standard argument (cf Lemma 3.3, page 15 of
\cite{B}) to again conclude that the CAT(0)--inequality holds.
\end{proof}

Next we consider some examples which show that in special cases both
questions can have positive answers even when the curvature is not
constant.

 \begin{Example}[\cite{GT}] Using the arithmetic
theory of quadratic forms, one can construct pairs $(V,W)$ as in
Example 4, where $V$ is a hyperbolic $n$--manifold. Consider the
quadratic form, $-\sqrt{2}x_0^2 + x_1^2 + \cdots + x_n^2$. Let
$\Gamma_n$ denote the group of isometries of this form with
coefficients in the ring of integers of ${\mathbb Q}(\sqrt{2})$. Then
$\Gamma_n$ is a discrete cocompact subgroup of $O(n,1)$. If $\Gamma $
is any torsion-free subgroup of finite index in $\Gamma_n$, then
$V={\mathbb H}^n/\Gamma $ is a closed hyperbolic $n$--manifold. Note
that for each $i$, $1\leq i\leq n$, the hyperbolic reflection $r_i$,
defined by $(x_0, \ldots , x_i, \ldots , x_n) \to (x_0, \ldots , -x_i,
\ldots , x_n) $ lies in $\Gamma_n$. If $\Gamma $ is required to be
normal in $\Gamma_n$, then $r_i$ descends to an isometric involution
$\overline{r_i}$ of $V$. Its fixed point set $V^{\overline{r_i}} $ is
then a locally convex, codimension-one submanifold of
$V$. Consequently, $W=V^{\overline{r_1}}\cap V^{\overline{r_2}}$ is a
locally convex, codimension-two submanifold of $V$.  If we further
require $\Gamma $ to lie in the identity component of $O(n,1)$, then
both $V$ and $W$ will be orientable. Finally, by passing to a subgroup
of index two if necessary, we may suppose, that each
$V^{\overline{r_i}} $ separates $V$ into two components. This implies
that $0=[W]\in H_{n-2}(V)$. Hence, we can obtain a pair $(V,W)$
satisfying all the conditions of Theorem A. It follows that the
induced (singular) metric on $V_k$ is nonpositively curved. This
metric is naturally piecewise hyperbolic. In \cite{GT} Gromov and
Thurston show that this metric on $V_k$ can be smoothed to a
Riemannian metric with pinched negative curvature. Moreover, if $n\geq
4$, $V_k$ is not diffeomorphic to a hyperbolic $n$--manifold for
infinitely many $k$.

The construction of Example 3 can be used to show that each
Gromov--Thurston example also admits nonpositively curved, piecewise
Euclidean metric.  Simply choose the net $S$ in $V$ so that its
inverse image $\widetilde{S}$ in ${\mathbb H}^n$ is
$\Gamma_n$--stable. The convex hull construction then gives a piecewise
Euclidean metric on $V$ with the property that each $\overline{r_i}$
acts isometrically. Hence, the fixed point sets $V^{\overline{r_1}} $
and $V^{\overline{r_2}} $ as well as their intersection $W$ will be
locally convex. The induced length metric on $V_k$ is piecewise
Euclidean and by Theorem A, it is nonpositively curved.
\end{Example}

 \begin{Example}[{\cite[page 383]{DJ}}]  In each dimension $\geq 5$ one
can construct a nonpositively curved, piecewise Euclidean metric on a
smoothable manifold $M^n$ such that its universal cover $X$ (a
Hadamard space) is not simply connected at infinity and hence, not
homeomorphic to ${\mathbb R}^n$.  On the other hand, by the
Cartan--Hadamard Theorem, if $M^n$ admitted a nonpositively curved
Riemannian metric, then its universal cover would be homeomorphic to
${\mathbb R}^n$. Thus, $M^n$ admits no such metric.
 \end{Example}

\section{Eberlein's Rigidity Theorem in the singular case}
Throughout this section, we shall make the following assumptions.

{\em Let $(X^{\ast }, d^{\ast })$ be a reducible global 
Riemannian symmetric space
with no compact or Euclidean factors
and of rank $r\geq 2$, and let $\Gamma^{\ast } $ be an irreducible
torsion-free uniform lattice on $X^{\ast }$.  Set $V^{\ast }=X^{\ast }/\Gamma^{\ast }$.}

{\em Let $(X,d)$ be a Hadamard space 
and let $\Gamma $  be a discrete
group of isometries on $X$ such that $\Gamma $ and $\Gamma^{\ast } $ are isomorphic as
groups and such that $V=X/\Gamma $ is  compact. We
further assume that $V$ is a topological manifold.}

The goal of this section is to show that the rigidity holds
even
without the smoothness assumption, namely, one has the following.

 \begin{TB} $(X,d)$ is isometric to $(X^{\ast },d^{\ast \ast
})$ for some $d^{\ast \ast }$ which differs from $d^{\ast }$ only by
constant multiples on the factors of $X^{\ast }$.
\end{TB}

By the singular analogue (\cite{BH}, Chapter IV, Theorem 6.1) of the
Schroeder's Splitting Theorem, we know that the theorem holds true if
$V^{\ast }$ is (the finite undercover of) the product of several such
irreducible compact locally symmetric spaces.

The proof will be a modification of the one given by Eberlein (\cite{E1},
\cite{E2}, \cite{CE}).

The key to understanding the geometry of a symmetric space $X^{\ast }$
(reducible or not) is the study of its maximal flats. Given a maximal
flat $F^{\ast }$ in $X^{\ast }$ and a point $x^{\ast } \in F^{\ast }$,
one considers the other maximal flats which contain $x^{\ast }$.  It
turns out that $F^{\ast }$ intersects these other maximal flats in a
union of hyperplanes.  These hyperplanes cut $F^{\ast }$ into
simplicial cones, any one of which is called a {\em Weyl chamber.}
Given a Weyl chamber and some other point $y^{\ast } \in X^{\ast }$,
there is a unique maximal flat $F^{\ast }$ which contains $y^{\ast }$
and a Weyl chamber parallel to the given one.

A geodesic line in $X^{\ast }$ is {\em singular\/} if it is contained in
more than one maximal flat.  It is {\em maximally singular\/} if it is
the intersection of the maximal flats which contain it, in other
words, if it is an extremal line in a decomposition of a maximal flat
into simplicial cones.

\subsection{The image flats}

In this subsection, we don't need the assumption that $X^{\ast }$ is
reducible.

First of all, choose a homotopy equivalence $\overline{f}\co V^{\ast
}\to V$ and an inverse $\overline{g}\co V\to V^{\ast }$.  Let $f$ and
$g$ be their lifts to the universal covers. Then $f$ and $g$ are
$\Gamma $--equivariant, uniformly continuous, and they are $(k,b)$
quasi-isometries, that is,
\[ \frac{1}{k} d^{\ast }(x^{\ast },y^{\ast }) \leq d(f(x^{\ast }), f(y^{\ast })) \leq k d^{\ast }(x^{\ast },y^{\ast }) \
\]
for any $x^{\ast }$, $y^{\ast }$ in $X^{\ast }$ with $d^{\ast
}(x^{\ast }, y^{\ast })\geq b$, and similarly for $g$.  We may also
choose a large constant $A$ such that for any $x^{\ast }$ and $y^{\ast
}$ in $X^{\ast }$,
\[ d^{\ast }(x^{\ast }, gf(x^{\ast })) \leq A, \mbox{ and } d(f(x^{\ast }),f(y^{\ast })) \leq
kd^{\ast }(x^{\ast },y^{\ast }) + A \] and similarly for any $x$ and
$y$ in $X$,
\[ d(x, fg(x)) \leq A, \mbox{ and } d^{\ast }(f(x),f(y)) \leq kd(x,y)
+ A.  \]

The group $\Gamma $ acts on $X^{\ast }$ via the isomorphism $\Gamma
\cong \Gamma^{\ast }$.  An $r$--flat $F^{\ast }\subseteq X^{\ast }$ is
called {\em $\Gamma$--compact\/} if $\Gamma F^{\ast } / \Gamma $ is
compact.  This is equivalent to the condition that $\Gamma$ contains
an abelian subgroup $L$ of rank $r$ such that $F^{\ast }$ is contained
in the minimum set of $L$ (which is defined as the intersection of the
minimum sets of all $\gamma^{\ast }$ in $L$).

Since $V$ is compact, all elements of $\Gamma $ are semisimple on
$X$. So $L$ is an abelian group of semisimple isometries on $X$,
therefore in the space $X$, the minimum set $MIN(L)={\mathbb R}^r
\times \Sigma \subseteq X$ splits isometrically, with $L$ acting as
identity on the convex set $\Sigma $ and as a lattice of translations
on ${\mathbb R}^r$.  Pick any $p\in \Sigma $ and let $F={\mathbb R}^r
\times \{ p\} $, then it is clear that both $ Hd(F, f(F^{\ast })) $
and $Hd(F^{\ast }, g(F))$ are finite, since all of the subsets
involved are $L$--invariant and have compact quotients. Therefore, we
have proved the following:

{\em For any $\Gamma $--compact $r$--flat $F^{\ast }\subseteq X^{\ast
  }$, there exists $r$--flat $F\subseteq X$ such that $Hd(F, fF^{\ast
  }) < \infty $.}

Note that $Hd(gF, F^{\ast })< \infty $, and the $r$--flat $F^{\ast }$
in $X^{\ast }$, which satisfies this is unique.  So, by the Quasi-flat
Theorem of Kleiner--Leeb \cite{KL} or Eskin--Farb (Corollary 7.4,
\cite{EF}), we know that there exists a constant $R'>0$, which is
independent of the choice of flats, such that $Hd(gF,F^{\ast }) \leq
R'$. By letting $R=kR'+2A$, we have the following lemma.

 \begin{Lemma} There exists a constant $R>0$ such that for any
   $r$--flat $F^{\ast }\subseteq X^{\ast }$, there is a $r$--flat
   $F\subseteq X$ with Hausdorff distance $Hd(F, fF^{\ast }) \leq R$.

   Furthermore, if $F_1$ is any $r$--flat in $X$ with $Hd(F_1, fF^{\ast
}) < \infty $, then \linebreak $Hd(F_1, fF^{\ast }) \leq R$. \end{Lemma}

 \begin{proof} From the above discussion, the statement holds when
$F^{\ast }$ is $\Gamma $--compact. Since the set of $\Gamma $--compact
$r$--flats is dense in the set of all $r$--flats in $X^{\ast }$, and any
sequence of $r$--flats in $X$ which all meet a compact set will have a
convergent subsequence, the general case follows.  \end{proof}

 We will call such a flat $F$ an {\em image flat\/} of $F^{\ast }$, and
 $F^{\ast }$ the preimage flat of $F$. Note that the preimage flat of
 $F$ is unique, while union of all the image flats of $F^{\ast }$
 forms a closed convex set which splits as a product $F\times \Sigma
 $, with $\Sigma $ convex and having diameter $\leq 2R$.

\subsection{The reducibility of $X$}

In this subsection, we want to show that the reducibility result of
\cite{E1} holds in the singular case as well. Most of the proof from
\cite{E1} goes through, with the following modifications.

First, we need to establish the so-called ``duality condition'' of
$\Gamma $ on $X$. Recall that $\Gamma $ satisfies the {\em duality
condition\/} on $X$, if for any geodesic $c$ in $X$, there exists a
sequence $\{\gamma_n\} \subseteq \Gamma $ such that $\gamma_n(x) \to
c(+\infty )$ and $\gamma_n^{-1}(x) \to c(-\infty )$. To establish this
we need the following lemma.

 \begin{Lemma} Let $c$ be a geodesic in $X$. Then there exists an
image flat $F$ of a flat $F^{\ast }\subseteq X^{\ast }$, such that
$d(c(t),F)$ is constant.
\end{Lemma}

\begin{proof}  Since there are $r$--flats in $X^{\ast }$ passing through any two given
points, we may choose $F_n^{\ast }\subseteq X^{\ast }$ passing through
$g(c(n))$ and $g(c(-n))$.  Let $F_n$ be an image flat of $F_n^{\ast
}$, $n=1,2,\ldots$. Since
\[ d(c(\pm n), F_n) \leq A+ d(fg(c(\pm n)), f(F_n^{\ast })) + R = A+R \
\]
by the convexity of the metric, we have $d(c(0),F_n)\leq A+R$, so by
passing to a subsequence, we may assume that $F_{n_k}\to F$ and $F$ is
an image flat of $F^{\ast }$, the limit of the $F_{n_k}^{\ast }$. We
have $d(c(\pm n_k), F)\leq A+R$ is bounded for all $k$; hence,
$d(c(t),F)$ is constant, by the convexity of the metric.
\end{proof}

\begin{Lemma}
  Any finite index subgroup of $\Gamma $ satisfies the duality
  condition on $X$.
\end{Lemma}

\begin{proof}  Fix a geodesic $c$ in $X$. By Lemma 2, there exists a flat $F^{\ast }$ and
an image flat $F$ such that $c$ has constant distance from $F$.  Let
$F_n^{\ast }$ be a sequence of $\Gamma $--compact flats which converges
to $F^{\ast }$. For each n, let $F_n$ be an image flat of $F_n^{\ast
}$, and denote by $L_n\subseteq \Gamma $ the corresponding lattice on
$F_n^{\ast }$ or $F_n$. Passing to subsequence if necessary, we may
assume that $F_n\to \widetilde{F}$, which is also an image flat of
$F^{\ast }$, hence $c$ has constant distance to $\widetilde{F}$.  Let
$\tilde{c}$ be the projection of $c$ onto $\widetilde{F}$. We can
choose a sequence of $L_n$--compact geodesics $c_n\subseteq F_n$ such
that $c_n \to \tilde{c}$. That is, there exists a sequence $\{
\gamma_k^{(n)} \} \subseteq L_n \subseteq \Gamma $ such that for each
$n$, $\gamma_k^{(n)}(c_n(0)) \to c_n(+\infty )$ and $\
(\gamma_k^{(n)})^{-1} (c_n(0)) \to c_n(-\infty )$ when $k\to \infty
$. By taking a diagonal subsequence in $\{ \gamma_k^{(n)}\} $, we know
that there exists sequence $\{ \phi_n\} \subseteq \Gamma $ such that $
\phi_n(x) \to c(+\infty )$ and $\phi_n^{-1}(x) \to c(-\infty )$ as
$n\to \infty $. \end{proof}

\begin{Lemma}[cf Lemma 2.4a  of \cite{CE}]\label{CE2.4a}
  Let $x \in X(\infty)$ be arbitrary and let $y \in X(\infty)$ be a
point that can be joined by a geodesic to $x$. If $z \in X(\infty)$ is
any point that can be joined to $x$ then $z \in \overline{\Gamma(y)}$.
\end{Lemma}

\begin{proof} With the duality condition established as above, the proof of
Lemma 2.4a in \cite{CE} can be easily modified to cover the singular
case. All we have to do is to estimate $\measuredangle_q(\phi_ny,
\phi_np)$ the same way for a point $p$ on $\gamma $ and allow $q$ to
be an arbitrary but fixed point in $X$, then $\measuredangle_q(w,z)
=0$ for any $q$ would imply $w=z$.  Here $w$ is any accumulation point
of $\phi_ny$.
\end{proof}

With the duality condition established, let us walk through the proofs
in sections 4 and 5 of \cite{E1}. Our metric space $(X, d)$ is no
longer smooth and we don't have the uniqueness of continuation of
geodesic segments.  Our starting point is the following proposition.

\begin{Prop}[cf Proposition 4.1 of \cite{E1}]\label{E4.1}
Let $\gamma^{\ast }$ be a maximally singular\break geo\-desic of $X^{\ast }$
and let $ F^{\ast }$ be an $r$--flat containing $\gamma^{\ast }$. Let
$F$ be an $r$--flat in $X$ such that $ Hd(f(F^{\ast }), F) < \infty
$. Then
\begin{enumerate}
\item the points $x=\lim_{t \to +\infty}(f \circ \gamma^{\ast })(t)$
and $y=\lim_{t \to -\infty}(f \circ \gamma^{\ast })(t)$ exist and are
distinct in $X(\infty)$
\item for any point $p$ in $F$ it follows that $
\measuredangle_p(x,y)=\pi $ and there exists the unique geodesic
$\gamma$ joining $x$ and $y$ passing through $p$, and $\gamma $ is
contained in $F$.
\end{enumerate}
\end{Prop}

\begin{proof} For Proposition~\ref{E4.1}, the proof given in the appendix of
\cite{E1} works.  The only place where care needs to be taken is with
regard to the paragraph on page 73, where the angles are used. In a
general Hadamard space, the angle function $\measuredangle_x(y,z)$,
while it is still continuous in $y$ and $z$, is only upper
semicontinuous with respect to the vertex $x$. Also, two different
geodesics from $x$ could form a zero angle (in fact they could even
share a segment).

Note that given two distinct points $u$, $v$ in $X( \infty )$, there
always exists a point $q\in X$ such that $\measuredangle _q(u,v)
>0$. Also, $\measuredangle _q(u,v) $ is continuous in $u$ or $v$ on
$\overline{X}=X \cup X(\infty )$, and satisfies the triangle
inequality
\[ \measuredangle _q(u,v) \leq \measuredangle _q(u,w) + \measuredangle
_q(w,v). \] With these properties, the paragraph on page 73 of
\cite{E1} can easily be modified and the proof is valid. \end{proof}

Now we consider the situation where $X^{\ast }=X_1^{\ast } \times
X^{\ast }_2$. Define $S^{\ast } \subset X^{\ast }(\infty)$ to be the
set of endpoints at infinity of the maximally singular geodesics of
$X^{\ast }$ and let $S_i^{\ast }=S^{\ast } \cap X_i^{\ast }$ for
$i=1,2$. ($S^{\ast }$ is the set of {\em maximally singular points\/} in
$X^{\ast }(\infty)$.)

By Proposition 1 the map $f$ defines an equivariant map
$\overline{f}\co S^{\ast } \to X(\infty)$. Let $S_i$ denote the
closure of $\overline{f}(S_i^{\ast })$ in $X(\infty)$.

\begin{Prop}[cf Proposition 4.3 of \cite{E1}]\label{E4.3}
Let $S_1$, $S_2$ be as above. Then:
\begin{enumerate}
  \item Each of the sets $S_1$, $S_2$ is invariant under $\Gamma$.
  \item A point $x \in S_i$ can only be joined by a geodesic to some
  other point in $S_i$.
\end{enumerate}
\end{Prop}

\begin{proof} (1) is obvious. (2) follows by taking closures from 
the fact that, 
for $x \in \overline{f}(S_i^{\ast })$, the point $y$ of
Proposition~\ref{E4.1} belongs to $S_i$, together with part (1) and
Lemma~\ref{CE2.4a}. \end{proof}

Now let us continue with the proofs in Appendix 2 of \cite{E1}. We
start with Lemma 5.1. Lemma 5.1 is needed in the proof of Proposition
5.1 and Sublemma 5.4b. We will state the following weaker version,
which is sufficient for our purpose.

\begin{Lemma}[cf Lemma 5.1 of \cite{E1}]\label{E5.1}
Suppose $\gamma (t)$ is a geodesic in the CAT(0) space $X$ and write
  $x=\gamma (\infty )$, $y= \gamma (-\infty )$, $p=\gamma (0)$. Denote
  by $L(p,x)$ and $B(p,x)$ the horosphere and horoball, respectively,
  which are centered at $x$ and pass through the point $p$.  Let
  $C=C_{xy}=\{ q\in X \mid \measuredangle _q(x,y) = \pi \} $.  Then
  $L(p,x)\cap L(p,y) = B(p,x)\cap B(p,y)$, $C$ is the union of all the
  geodesics parallel to $\gamma $, $C$ is a closed convex subset of
  $X$, and $C$ splits as $C=\gamma \times C'$ for a closed convex subset
  $C' \subseteq L(p,x)\cap L(p,y)$.
\end{Lemma}

\begin{proof} In the singular case, it is still true that if two geodesic rays
meet at $q$ with angle $\pi $, then they form a geodesic. Also, two
parallel geodesics bound a flat strip (the `Flat Strip Lemma' or
`Sandwich Lemma', see for example \cite{BGS} or \cite{BH}).  By the
definition of Busemann function (see \cite{BGS}) and the Law of Cosine
in a comparison triangle, it is not hard to see that if $q\in B(p,x)$,
then $\measuredangle _p(x,q) \leq \frac{\pi }{2}$. Also, if $q\in
B(p,x) \setminus L(p,x)$, then $ \measuredangle _p(x,q) < \frac{\pi
}{2}$. From this it follows that $B(p,x)\cap B(p,y) \subseteq
L(p,x)\cap L(p,y)$. \end{proof}

For each point $p \in X$, let $B_p$ denote the smallest convex subset
of $X$ which contains all the maximally singular geodesics through $p$
with endpoints in $S_1$.

\begin{Prop}\label{Ep5.1}
  With notation as in Lemma~\ref{E5.1}, suppose $x\in S_2$. Then $B_p$
  is contained in an $R'$--neighborhood of $C_{xy}$, for some constant
  $R'$.
\end{Prop}

\begin{proof} This is essentially the proof of Proposition 5.1 of \cite{E1}, but
since we do not assume that $C_{xy}$ is a proper subset of $X$, we get
a weaker result.
\end{proof}

\begin{Prop}[cf Proposition 5.2 of \cite{E1}]\label{E5.2}
  If $r\in B_p$ then $B_r \subseteq B_p$.
\end{Prop}

\begin{proof} The proof of Proposition 5.2 of \cite{E1}  is valid in the
 singular case as well.
\end{proof}

Before stating the next proposition we need the following four lemmas.

\begin{Lemma}[cf Lemma 5.5a of \cite{E1}]\label{E5.5a}
  Let $p$, $r$ be any two points of $X$. Then $d(q,B_r)\leq d(p,B_r)$
  for every point $q \in B_p$. Moreover suppose that $p \not\in B_r$
  and let $P\co X \to B_r$ denote the nearest point projection. Then
  for every $x \in S_1$ and any geodesic $\gamma$, with $\gamma(0)=p$
  and $\gamma(\infty)=x$, the projection $ (P\circ\gamma)(t) $ is a
  unit speed geodesic and $\gamma$ and $P\circ\gamma$ bound a flat
  strip in X.
\end{Lemma}

\begin{Lemma}[{\cite[page 25, Corollary 5.8 (i)]{B}}]\label{E5.4a}
   Let $p,\ q \in X$ and $x \in X(\infty )$. If $
\measuredangle_p(x,q) + \measuredangle_q(x,p) = \pi $, then the three
points span a flat strip.
\end{Lemma}

\begin{Lemma}[cf Sublemma 5.4b of \cite{E1}]\label{E5.4b}
Let $\gamma $ be a geodesic of $X$ such that 
$\gamma [0,+\infty)\subseteq B_p$, 
where $ p=\gamma (0)$. Let $ z = \gamma (+\infty)$. Then there exists $x \in S_1$ and $t \in {\mathbb R}$ such
that $\measuredangle_{\gamma(t)}(x,z)\not=\pi/2$.
\end{Lemma}

\begin{proof} The proof of Sublemma 5.4b of \cite{E1} holds  true using
Lemma~\ref{E5.1} and Lemma~\ref{E5.4a}, instead of Lemma~5.1 and
Sublemma~5.4a of \cite{E1}.
\end{proof}

\begin{Lemma}[cf Lemma 5.4 of \cite{E1}]\label{EL5.4}
For every point $p \in X$ there exist a number $A>0$ such that
$d(p,B_r)<A$ for every $r \in B_p$.
\end{Lemma}

\begin{proof} Again, the proof of  Lemma 5.4 of \cite{E1} works using
Lemma~\ref{E5.5a} instead of Lemma~5.5a of \cite{E1} and the fact that
the limit of a convergent sequence of (half) flats is still a (half)
flat.
\end{proof}

\begin{Prop}[cf Proposition 5.4 of \cite{E1}]\label{E5.4}
  There exists a point $p\in X$ such that $B_r = B_p$ for every point
$r\in B_p$.
\end{Prop}

\begin{Prop}[cf Proposition 5.5 of \cite{E1}]\label{E5.5}
 For every point $p\in X$ and every point $r\in B_p$, $B_r = B_p$.
\end{Prop}

\begin{Prop}[cf Proposition 5.6 of \cite{E1}]\label{E5.6}
 For every point $p\in X$ the set $B_p$ is complete, totally geodesic
 subspace of $X$.
\end{Prop}

\begin{proof} The proof of Proposition 5.6 of \cite{E1}  goes through without
change, except for the last paragraph of the proof on page 63, where
angles are again used. What is needed to complete the proof is the
following lemma.

\begin{Lemma}\label{5.6}  Consider a geodesic triangle with vertices
$p$, $q$, and $r$. Let $\sigma (t)$, $t\in [0,
  d(p, q)]$, be the geodesic segment from $p$ to $q$, with $\sigma
  (0)=p$. If there is a sequence $\{ t_i\} $ of small positive numbers
  that converges to $0$, such that $d(r,p) \leq d(r, \sigma (t_i) $
  for each $i$, then $\measuredangle_{p} (q, r) \geq \frac{\pi }{2}$.
\end{Lemma}

\begin{proof} Let $\gamma (t)$ be the geodesic segment from $p$ to $r$, with
$\gamma (0)=p$, $\gamma (l)=r$, where $l=d(p,r)$. Denote by $d_{st}$
the distance between $\gamma (s)$ and $\sigma (t)$. Then the function
\[ f(t,s) = (t^2+s^2 - d^2_{st})/2st  \]
is monotonically decreasing in $s>0$ and $t>0$, and it's limit when
$s\to 0$, $t\to 0$ is the cosine of the angle $\measuredangle_{p} (q,
r)$.

Now assume that $\measuredangle_{p} (q, r) < \frac{\pi }{2}$. Then
there exist $ t >0$ and $s>0$ such that $f(t,s) \geq \epsilon >0$. By
the monotonicity of $f$ in $t$, we have
\[ (t'^2 + s^2 - d^2_{st'})/ 2st' \geq \epsilon  \
\]
for any $0<t'\leq t$. So for $t'$ sufficiently small, we have $s >
d_{st'}$.  Hence for $t_i$ sufficiently small,
\[ d(r,p) = d(r, \gamma (s)) +s > d(r, \gamma (s)) + d_{st_i} \geq d(r,
\sigma (t_i) \] which contradicts our assumption. \end{proof}

To complete the proof of Proposition 7, we apply the lemma to the
triangle with vertices $q_1$, $\overline{q}$ and $r$.  We have that
two of the internal angles are at least $\frac{\pi }{2}$. So the sum
of the internal angles is at least $\pi $, which implies that it must
be a triangle in ${\mathbb R}^2$, with two right internal angles. This
is clearly impossible. \end{proof}

Summarizing the above discussions, we can now state the Eberlein's
result in \cite{E1} in the singular case as the following:

 \begin{TC} Let $X^{\ast }=X^{\ast }_1\times \cdots \times
  X^{\ast }_k$ be the decomposition into irreducible factors, $k\geq
  2$. For $1\leq i\leq k$, denote by $S_i^{\ast } \subseteq X_i^{\ast
  }(\infty )$ the set at infinity of the maximally singular geodesics in
  $X_i^{\ast }$. Then the following holds:

  {\rm(1)}\qua For any maximally singular geodesic $c^{\ast }\subseteq X^{\ast
  }$, there exists a geodesic $c\subseteq X$ such that the Hausdorff
  distance $Hd(c, f(c^{\ast }))\leq R_1$ for some uniform constant
  $R_1$.

  {\rm(2)}\qua If $c^{\ast }$ is the intersection of $r$--flats $F_1^{\ast },
\ldots , F_m^{\ast }$, and $c_1$ is a geodesic parallel to all the
image flats $F_1, \ldots , F_m$, then $c_1$ is parallel to $c$.

  {\rm(3)}\qua Let $S_i\subseteq X(\infty )$ be the image of $S_i^{\ast }$
guaranteed by (1). For any $x\in X$, denote by $A_i(x)$ the union of
all the geodesic rays from $x$ to the points in $S_i$, and let
$B_i(x)$ be the smallest closed convex subset of $X$ that contains
$A_i(x)$. Then each $B_i(x)$ is totally geodesic in $X$, and $ B_i(x)$
is parallel to $B_i(y)$ for any $x,y\in X$.

  {\rm(4)}\qua For any $x\in X$ and $z\in S_j$, let $c$ be a geodesic with
$c(0)=x$ and $c(\infty )=z$. Then for $i\neq j$, $B_i(x) \subseteq
N_{R'}(P_c)$, the $R'$--neighborhood of $P_c$. Here $P_c$ denotes the
set of points $y\in X$ where there is a geodesic through $y$ that is
parallel to $c$. 
\end{TC}

\subsection{Transitivity of the isometry group  on $X$}
Our next goal is to establish the fact that the isometry group $I(X)$
of $X$ acts transitively on $X$, hence $I(X)$ is a Lie group by the
theorem of Montgomery and Zippin \cite{MZ}.

Let $X^{\ast }=X_1^{\ast } \times \cdots \times X_k^{\ast }$ be the
decomposition into irreducible factors. By assumption, $k\geq 2$.  In
order to distinguish the actions of $\Gamma $ on $X^{\ast }$ and $X$,
we will denote the lattice on $X^{\ast }$ by $\Gamma^{\ast }$.  By
passing to a sublattice if necessary, we may assume that the
irreducible lattice $\Gamma^{\ast } \subseteq I_0(X^{\ast
})=I_0(X_1^{\ast }) \times \cdots \times I_0(X_k^{\ast })$, where
$I_0$ stands for the identity component of the isometry group. So each
$\gamma^{\ast }\in \Gamma^{\ast }$ is in the form $\gamma^{\ast }= (
\gamma_1^{\ast }, \ldots , \gamma_k^{\ast })$. Denote by
$\Gamma_i^{\ast }$ the projection of $\Gamma^{\ast }$ on $X_i^{\ast
}$. Since $\Gamma^{\ast }$ is assumed to be irreducible, by Borel
Density Theorem, we know that for each $1\leq i\leq k$, $\Gamma^{\ast
} \to \Gamma_i^{\ast }$ is injective, and $\Gamma_i^{\ast }$ is dense
in $I_0(X_i^{\ast })$.

From the discussion in the last subsection, we know that for each
$1\leq i\leq k$, $X$ is foliated by parallel, totally geodesic leaves
$B_i$. By the Sandwich Lemma, for each $i$, $X$ splits isometrically
as $B_i \times Y_i$.  Basically, we want to show that these foliations
are mutually perpendicular, and hence, that $X$ has a product
structure corresponding to that of $X^{\ast }$. But we can only prove
this in the orthogonal complement of the Euclidean factor of $X$.

Denote by $C(X)$ the set of all Clifford translations in $I(X)$, the
isometry group of $X$.  It is a normal abelian subgroup of $I(X)$.
$X$ splits isometrically $X={\mathbb R}^s \times Y$, $s\geq 0$, such
that each $\phi \in C(X)$ has the form $(T, id)$ under this splitting
with $T$ a translation, and so that $Y$ does not contain any Euclidean
factor. Furthermore, any isometry of $X$ respects this splitting:
$I(X)=I({\mathbb R}^s)\times I(Y)$. (cf (6.2) of \cite{BGS}).

Our first goal is to show that, for each $1\leq i\leq k$, the leaves
of $B_i$ all respect the product structure $X={\mathbb R}^s\times
Y$. To be more precise, for any $x=(x_0,y) \in X$, $B_i(x) =
B_{i0}(x_0) \times B_{iY}(y)$, and $B_{i0}$ (or $B_{iY}$) is a
parallel totally geodesic foliation in ${\mathbb R}^s$ (or $Y$).

As before, for each $1\leq i\leq k$, let $S_i^{\ast }\subseteq
X_i^{\ast }(\infty )$ denote the set at infinity of the maximal
singular geodesics in $X_i^{\ast }$ ($1\leq i\leq k$), and
$S_i\subseteq X(\infty )$ the image guaranteed by part (1) of Theorem
C. Again denote by $A_i(x)$ the union of all geodesic rays from $x$ to
the points in $S_i$, and $B_i(x)$ the smallest closed convex set
containing $A_i(x)$.

Now consider an image flat $F\subseteq X$. Let us denote by $F_0$ and
$F_Y$ the projection of $F$ in ${\mathbb R}^s$ and $Y$,
respectively. If we fix a point $x=(x_0,y) \in F$ and regard $F_0
\subseteq {\mathbb R}^s \times \{ y\} $, $F_Y \subseteq \{ x_0 \}
\times Y$, then we have the following lemma.

\begin{Lemma}
 For any image flat $F$ in $X$, $F=F_0 \times F_Y$.
\end{Lemma}

\begin{proof} Since any image flat is parallel to an image flat that is the
limit of a sequence of $\Gamma $--compact image flats, it suffices to
prove the lemma for $\Gamma $--compact image flats. Let $L\subseteq
\Gamma $ be a rank $r$ abelian subgroup which acts as a lattice on
$F$. For any $\gamma \in L$, write $\gamma = (\gamma_0, \gamma_Y)$
under the splitting $X={\mathbb R}^s \times Y$. Then $\gamma_0$ is a
translation on $F_0$. This is true because, for any geodesic $c_0$ in
$F_0$, $c_0$ is the projection of some geodesic $c$ in $F$. Since
$\gamma $ maps $c$ to a parallel geodesic in $F$, $\gamma_0(c_0)$ is
parallel to $c_0$. Let $T$ be a translation in $F_0$.  As $T$ commutes
with $\gamma_0$ for all $\gamma \in L$, $TF$ is also invariant under
$L$. Since we have uniform bound in Hausdorff distance between $g(F)$
or $g(TF)$ and $F^{\ast }$, the preimage flat, we know that the
Hausdorff distance between $F$ and $TF$ must be bounded by a uniform
constant. Therefore, $F_0$ must be contained in $F$, since otherwise
we can choose a translation in $F_0$ that sends $F$ to a parallel $TF$
with arbitrarily large $Hd(F,TF)$.
\end{proof}

As a corollary of this, we have the following.

\begin{Lemma} For each $1\leq i\leq k$, the set $S_i \subseteq
{\mathbb R}^s(\infty )\cup Y(\infty )$. So the foliation $B_i$ is a
    product foliation $B_{i0} \times B_{iY}$. Furthermore, $Y=B_{1Y}
    \times \cdots \times B_{kY}$.
\end{Lemma}

\begin{proof}  Since any maximally singular geodesic $c^{\ast } 
\subseteq X_i^{\ast }$ is an intersection of $r$--flats,  there
exist $r$--flats $F_1^{\ast }, \ldots , F_m^{\ast }$ with  $c^{\ast
}=F_1^{\ast }\cap \ldots \cap F_m^{\ast }$.  Let $F_j$ be an image
flat of $F_j^{\ast }$, $1\leq j\leq m$. Then part (1) of Theorem C
says that there exists a geodesic $c\subseteq X$ which is parallel to
all $F_j$, and any geodesic $\tilde{c}$ parallel to all $F_j$ must be
parallel to $c$. In the product space $X={\mathbb R}^s\times Y$, $c$
is parallel to $F_j$ if and only if $c_0$ is parallel to $F_0$ and
$c_Y$ is parallel to $F_Y$. Here the subscripts denote the
projections.  So if $c$ is not parallel to a leaf of ${\mathbb R}^s$
or $Y$, then any line in the plane $c_0\times c_Y$ would be parallel
to all those $F_j$, which is impossible. So, $c$ must be parallel to
either ${\mathbb R}^s$ or $Y$. That is, $S_i$ is contained in the
disjoint union of ${\mathbb R}^s(\infty ) $ and $ Y(\infty )$.

From this and the definition of $B_i$, it is clear that for each $i$,
$B_i = B_{i0} \times B_{iY}$ is a product foliation.

Next we claim that for $i\neq j$, $B_{iY}$ and $B_{jY}$ are
perpendicular.  Take any geodesic $c$ with $c(\infty )\in S_j \cap
Y(\infty )$. Then any leaf $B_i(x)$ is contained in the neighborhood
$N_{R'}(P_c)$ by part (4) of Theorem C. Here $P_c$ stands for the
union of all geodesics parallel to $c$.  So if the projection
$\overline{c}$ of $c$ onto $B_{iY}$ is not a point, any point in
$B_{iY}$ will be within distance $R'$ to another point in $B_{iY}$
where there is a geodesic parallel to $\overline{c}$ passing through.
That is, $B_{iY} \subseteq N_{R'}(P) $, where $P$ is the union of all
geodesics in $B_{iY}$ that is parallel to $\overline{c}$. Since
$B_{iY}$ is geodesically complete, and $P\subseteq B_{iY}$ is convex, we
know that $B_{iY}=P$. That is, $B_{iY}$ is foliated by geodesics
parallel to $\overline{c}$, so $B_{iY}$, hence $Y$ will contain a
Euclidean factor, a contradiction. So geodesics from $x\in B_{iY}$ to
any point in $S_j\cap Y(\infty )$ must be perpendicular to
$B_{iY}$. So $A_j(x)$ must be contained in $Z(x)$, the leaf through
$x$ of the orthogonal complement of $B_{iY}$ in $Y$, where $Y=B_{iY}
\times Z$ by the Sandwich lemma. Since $Z(x)$ is convex, we have
$B_{jY}\subseteq Z$, hence $B_{iY}$ is perpendicular to $B_{jY}$.  On
the other hand, since any geodesic in $X$ is parallel to an image
flat, which is contained the span of $B_1$ through $B_k$, so $Y$ is
spanned by those $B_{iY}$, $1\leq i\leq k$. \end{proof}

Note that if $X_i^{\ast }$ is of rank $1$, then $S_i^{\ast }$, hence
$S_i$ is connected, so either $S_i\subseteq {\mathbb R}^s(\infty )$ or
$S_i\subseteq Y(\infty )$. In this case $B_i$ is entirely in ${\mathbb
R}^s$ or $Y$.

Write $ X={\mathbb R}^s \times B_{1Y} \times \cdots \times B_{kY}$.
Each $\gamma \in \Gamma $ can be written in the form $\gamma =
(\gamma_0, \gamma_{1Y}, \ldots , \gamma_{kY}) = (\gamma_0, \gamma_Y)$,
since it preserves all the foliations $B_i$ and respects the product
structure $X={\mathbb R}^s \times Y$. Denote by $\Gamma_Y$ the group
of all such $ \gamma_Y $.  Set $\widetilde{\Gamma^{ \ast }} =
\Gamma_1^{\ast } \times \cdots \times \Gamma_k^{\ast } $, where
$\Gamma_i^{\ast }$ is the projection of $\Gamma^{\ast }$ onto the
irreducible factor $X_i^{\ast }$.  The group $ \widetilde{\Gamma^{\ast
}}$ is dense in $I_0(X^{\ast })$.

We want to construct a surjective, continuous map $\tilde{f}\co
X^{\ast }\to Y$ which is $(\widetilde{\Gamma^{\ast }},
\Gamma_Y)$--equivariant. Then $\Gamma_Y$ will have a dense orbit in
$Y$. Because $Y$ is locally compact, any sequence of isometries $\{
\phi_n \} $ on $Y$ such that $\phi_n(y)$ converges for some $y\in Y$
will have a subsequence that converges to an element in $I(Y)$.
Therefore, $I(Y)$ acts transitively on $Y$. Hence, it is a Lie group
by the theorem of Montgomery and Zippin \cite{MZ}. Thus, $I(X)$ is
also a Lie group.

If $Q$ is a bounded subset in a Hadamard space $Y$, we will denote by
$\tau_Q $ the infimum of the radii of the closed balls containing $Q$.
There exists an unique point $y\in Y$ such that $Q\subseteq
B_{\tau_Q}(y)$. This point is called the {\em circumcenter\/} of $Q$ and
is denoted by $\star (Q)$. We refer the reader to page 26 of \cite{B}
for more details.

\begin{Lemma} Suppose $P$ and $Q$ are two bounded subsets in a Hadamard
  space $Y$ with Hausdorff distance $Hd(P,Q)=h$. Then their
circumcenters satisfy $$ d(\star(P), \star(Q)) \leq \sqrt{
h(\tau_P+\tau_Q+h) }.$$
\end{Lemma}

\begin{proof} Write $x=\star(P)$ and $y=\star(Q)$. Then we have $P \subseteq
B_{\tau_P}(x)\cap B_{\tau_Q+h}(y)$. Let $z$ be the mid point of the
geodesic segment from $x$ to $y$. Then since $Y$ is Hadamard, any
$p\in P$ will satisfy
\[ d^2(p,z) \leq \frac{1}{2} d^2(p,x) + \frac{1}{2} d^2(p,y) -
\frac{1}{4} d^2(x,y) \leq \frac{1}{2} \tau_P^2 + \frac{1}{2}
(\tau_Q+h)^2 - \frac{1}{4} d^2(x,y). \] Since this is true for any
$p\in P$, the far right hand side of the above inequality must be
bigger than $\tau_P^2$, hence
\[ d^2(x,y) < 2\tau_Q^2 - 2\tau_P^2 + 4\tau_Q h + 2h^2. \]
Add this with the similar inequality obtained by reversing the role of
$P$ and $Q$, we get
\[ d^2(x,y) < 2h(\tau_P + \tau_Q +h).  \]
If one replace $z$ by the points on the geodesic segment from $x$ to
$y$ that are very closed to $x$ or $y$, then the coefficient $2$ in
the right hand side can be removed. \end{proof}

We are now ready to construct $\tilde{f}$. Let $\pi $ be the
projection map from $X$ onto $B_{1Y}$. For any $x^{\ast } \in
X_1^{\ast }$, write $Q^{\ast } (x^{\ast }) = \{ x^{\ast }\} \times
X_2^{\ast } \times \cdots \times X_k^{\ast }$, and $Q(x^{\ast }) = \pi
f(Q^{\ast }(x^{\ast }))$.  Define $\widetilde{f_1}(x^{\ast }) = \star
( Q(x^{\ast })) $ to be the circumcenter of $Q(x^{\ast })$. For this
we need to show that $Q(x^{\ast })$ is always a bounded subset of
$B_{1Y}$. Fix a point $p^{\ast }\in Q^{\ast }$. For any $q^{\ast }\in
Q^{\ast }$, there exists a $r$--flat $F^{\ast }\subseteq X^{\ast }$
containing both $p^{\ast }$ and $q^{\ast }$. So, if $r_1$ denotes the
rank of $X_1^{\ast }$, then there are $m\leq r-r_1$ maximally singular
geodesics in $F^{\ast }$, denoted by $l_1^{\ast }, \ldots , l_m^{\ast
}$, each is perpendicular to $X_1^{\ast }$, $p^{\ast }\in l_1^{\ast
}$, $q^{\ast }\in l_m^{\ast }$, and $l_i^{\ast }\cap l_{i+1}^{\ast }
\neq \phi $ for $1\leq i\leq m-1$. For each $l_i^{\ast }$, by Theorem
C, there exists geodesic $l_i$ in $X$ perpendicular to $B_{1Y}$, such
that the Hausdorff distance between $l_i$ and $f(l_i^{\ast })$ is
bounded by an uniform constant $R_1$. Therefore $d(\pi f(q^{\ast }),
\pi f(p^{\ast })) \leq 2m R_1 \leq 2rR_1$, and $\tau_Q \leq
2rR_1$. Since $f$ is uniformly continuous, and $Hd(Q(x^{\ast }),
Q(y^{\ast })) \leq d(x^{\ast },y^{\ast })$, so Lemma 6 implies that
$\widetilde{f_1}\co X_1^{\ast }\to B_{1Y}$ is continuous. It is also
clear that for any $\gamma \in \Gamma $, we have $\widetilde{f_1}
(\gamma_1^{\ast }(x^{\ast })) = \gamma_{1Y} (\widetilde{f_1}(x^{\ast
}))$. Now if we define $\widetilde{f_i}\co X_i^{\ast } \to B_{iY}$
similarly for each $1\leq i\leq k$ and then take the product, we get a
continuous map $\tilde{f}\co X^{\ast } \to Y$ which is $(\widetilde{\Gamma^{\ast }}, \Gamma_Y) $--equivariant, where $
\widetilde{\Gamma^{\ast }} = \Gamma_1^{\ast } \times \cdots \times
\Gamma_k^{\ast } $. This implies the following.

 \begin{Lemma} The map $\tilde{f}\co X^{\ast }\to Y $ is surjective.
Hence $\overline{\Gamma_Y}$ acts transitively on $Y$, and $I(X)$ is a
Lie group.
\end{Lemma}

Note that for any $x^{\ast }\in X^{\ast }$, we have
$d(\tilde{f}(x^{\ast }), \pi_Yf(x^{\ast })) \leq \sqrt{k}2rR_1$. So
for any $y\in Y$,
\[ d(y, \tilde{f}\circ g(0,y)) \leq A + 2r\sqrt{k}R_1. \]
Therefore $\tilde{f}$ is surjective by the following lemma.

\begin{Lemma} If a Hadamard space $Y$ is a topological manifold and if
$h\co Y\to Y$ is a continuous map within bounded distance from the
  identity map, then $h$ is surjective.
\end{Lemma}

\begin{proof} The geodesic segment from $h(y)$ to $y$ gives a proper homotopy
between $h$ and the identity map.  Hence, $h$ induces the identity map
on the top dimensional cohomology group with compact support,
$H^n_c(Y)$, $n=$dim$(Y)$.  The result follows.
\end{proof}

\subsection{Proof of the Rigidity Theorem}
Now we can apply the techniques in \cite{E2}, \cite{E3}, \cite{E4} and
\cite{CE} to finish the proof of Theorem B. From the previous
subsection, $I(X)$ is Lie group. This implies that $X$ has no
Euclidean factor, that is, we have the following lemma.

\begin{Lemma} $X=Y=B_1 \times \cdots \times B_k$.
\end{Lemma}

\begin{proof} First, notice that $\Gamma $ does not contain any Clifford
translations. This is because the set of Clifford translations in
$\Gamma $ forms a normal abelian subgroup. By the Main Theorem of
\cite{E3},
$\Gamma^{\ast }$ does not contain any nontrivial normal abelian
subgroup, since $X^{\ast }$ has no Euclidean de Rham factor.

Now let us assume that $X={\mathbb R}^s\times Y$ for $s>0$. We want to
derive a contradiction. Note that $Y$ is not a point, for otherwise
$\Gamma $ would contain a translation by the Bieberbach Theorem.

Denote by $\Gamma_{\mathbb R}$ and $\Gamma_Y$ the projection of
$\Gamma $ on ${\mathbb R}^s$ and $Y$, respectively. Since $I(X)$ is a
Lie group, the proof of Lemma A in \cite{E3} says that $\Gamma_Y$ is
discrete. Let $N$ be the kernel of $\Gamma \to \Gamma_Y$. The proof of
Theorem 4.1 of \cite{E4} says that either $\Gamma_{\mathbb R}$ is
discrete or $N$ contains a Clifford translation. Since $\Gamma^{\ast
}$ is assumed to be irreducible, $\Gamma_{\mathbb R}$ can not be
discrete. Since $\Gamma $ contain no Clifford translation, so does the
subgroup $N$. This completes the proof. (Note that in this argument,
Proposition 2.3 of \cite{CE} is used. But it is obviously valid in the
singular case.) \end{proof}

In our next argument we also need an extension of Theorem 2.4 of
\cite{CE} to singular spaces. The proof there is based on four
lemmas. We already established singular version of Lemma~2.4a --- see
Lemma~\ref{CE2.4a}.  The same sort of modification works for
Lemma~2.4c. The proof of Lemma~2.4b holds without change. For the
proof of Theorem~2.4 the following replacement of Lemma~2.4d suffices.

\begin{Lemma}
  Let $A \subseteq I(X)$ be a nontrivial abelian subgroup such that
  normalizer D of $A$ in $I(X)$ satisfies the duality condition. Then
  $I(X)$ contains a Clifford translation.
\end{Lemma}

\begin{proof} {\sl(cf proof of Lemma 2.4d of \cite{CE})}\qua
Let $L(A) \subseteq X(\infty)$
denote the limit set of $A$.  By Lemma~2.4b of \cite{CE}, $L(A)$ is
not empty.  Let $x\in L(A)$, then by Lemma~2.4c there exist a unique
point $y\in L(A)$ such that $x$ can be joined to $y$. If $z\in
H(\infty)$ is a point that can be joined to $y$ then by
Lemma~\ref{CE2.4a} $z\in \overline{D(x)} \subseteq L(A)$ and therefore
$z=x$.  So $X$ admits Clifford translations along the geodesics
joining $x$ to $y$ by Sandwich Lemma.
\end{proof}

Denote by $\Gamma_i = \Gamma_{iY}$ the projection of $\Gamma $ on
$B_i$, and write $G_i = (\overline{\Gamma_i})_0 \subseteq I_0(B_i)$
the identity component of its closure. By the proof of Lemma 3.1 of
\cite{E2}, which essentially uses Theorem 2.4 of \cite{CE} (cf Lemma
17), we know that each $G_i$ is a centerless semisimple Lie group of
noncompact type. Fix a maximal compact subgroup $K\subseteq G_i$. By
taking the circumcenter of an orbit, we know that $K$ will fix some
point in $B_i$.  Let $p\in B_i$ be a fixed point of $K$. Since $B_i$
is locally compact, the isotropy subgroup of $G_i$ at $p$ is compact,
and equals to $K$ since $K$ is maximal compact. So the transitivity of
the action of $G_i$ on $B_i$ gives a homeomorphism between $B_i$ and
$G_i/K$. The metric $d$ on $B_i$ is a $G_i$--invariant Hadamard
metric. By the following lemma, $d$ must be smooth, so $(B_i,d)$ is a
global Riemannian symmetric space of noncompact type.  Therefore $(X,
d)$ is symmetric, and isometric to $(X^{\ast },d^{\ast \ast })$ by the
theorem of Mostow \cite{M}. This completes the proof of Theorem B.

\begin{Lemma}\label{G-inv}
 Suppose $G$ is a semisimple Lie group of noncompact type and that $K$
  a maximal compact subgroup.  If $d$ is a Hadamard metric on $G/K$
  that is $G$--invariant, then $(G/K,d)$ is Riemannian symmetric space
  of noncompact type.
\end{Lemma}

\begin{proof} If $G$ is not simple, we may take a lattice $\Gamma
\subseteq G$ and consider $X^{\ast }=G/K$ equipped with $d^{\ast }$
symmetric and $X=G/K$ equipped with $d$. By the previous arguments, we
know that $(X,d)$ will be a product space of $(G_i/K_i,d_i)$, where
those $G_i$ are the simple factors of $G$. Therefore we may assume
that $G$ is simple.

First consider the case when the rank of $G$ is $1$. In this case, the
geodesic sphere in $d$ is a single $K$ orbit, so it coincides with a
geodesic sphere with the same center under $d^{\ast }$, the symmetric
metric. For any two points $p$, $q\in X$, let $m$ be the midpoint
under the metric $d$.  Then there are two geodesic spheres under
$d^{\ast }$, centered at $p$ and $q$, respectively, so that $m$ is
their unique intersection point. This implies that $m$ must be on the
geodesic segment from $p$ to $q$ with respect to the metric $d^{\ast
}$. So the geodesics in $d$ and $d^{\ast }$ have the same images.

Fix $p$ and $q$ with $d(p,q)=1$. Let $d^{\ast }$ be the symmetric
metric on $X=G/K$ such that $d^{\ast }(p,q)=1$. Let $\varphi_t$ be a
one parameter subgroup of $G$ such that $\gamma (t)= \varphi_t(p)$ is
the unit speed geodesic under $d^{\ast }$ with $\gamma (1)=q$.

By considering $\varphi_{\frac{1}{2^n}}$, we know that $d(\gamma (0),
\gamma (t))= d^{\ast }(\gamma (0), \gamma (t))$ for any rational
number $t$, whose denominator is a power of two.  Hence, it also holds
for any $t\in {\mathbb R}$ by continuity. For any $p', q'\in X$, and
for any $\varphi \in G$ such that $\varphi (p') = \gamma (0) $ and
such that $\varphi (q')$ lies in the image of $\gamma $, we have
$d(p',q') = d^{\ast }(p',q')$. So $d$ is symmetric in this case.

Now assume $G$ is simple with rank $r\geq 2$. We need the following.

\begin{Sublemma}  Suppose $G$ is a connected, centerless,
  semisimple Lie group of noncompact type, and $K$ a maximal compact
  subgroup. Let $\mathfrak{ g}=\mathfrak{ k} + \mathfrak{ p}$ be the
  Cartan decomposition of the Lie algebras. Suppose $\mathfrak{ a}$ is
  a maximal abelian subalgebra in $ \mathfrak{ p}$, and
  $A=\mbox{exp}(\mathfrak{ a})$.  Denote by $M=Z_K(A)$ the centralizer
  of $A$ in $K$, and denote by $Y$ the fixed point set in $G/K$ of
  $M$. Then $Y$ is a totally geodesic submanifold in $G/K$ containing
  the maximal flat $F=Ax_0$, where $x_0=K$. Denote by $r$ the rank of
  $G/K$.
\begin{enumerate} 
	\item $M=\{ 1 \} $ when and only when $G/K$
  is the product of $r$ copies of the hyperbolic plane ${\mathbb
  H}^2$.
	\item If $G$ is simple, then either $Y = ({\mathbb H}^2)^r$, or
$Y=F\cong {\mathbb R}^r$, depending on whether $G/K$ is Hermitian
symmetric or not.
\end{enumerate} 
\end{Sublemma}

Now let us use Sublemma to finish the proof of Lemma~\ref{G-inv}.
Let $F\subseteq B=G/K$ be a maximal flat under $d^{\ast }$, with $K$
the isotropy subgroup at $x_0\in F$. First we claim that $F$ is
$d$--convex.

Let $\mathfrak{g}= \mathfrak{k} + \mathfrak{p}$ be the Cartan
decomposition of the Lie algebras of $G$ and $K$, and $ \mathfrak{ h}$
a maximal abelian subalgebra contained in $ \mathfrak{p}$.  Write
$A=\mbox{exp} \mathfrak{ h}$ so that $F= Ax_0$. As in Sublemma,
denote by $M$ the centralizer of $A$ in $K$, and denote by $Y$ the fixed point
set of $M$ in $B$.  Then $Y\subseteq B$ is totally geodesic in $B$
under $d^{\ast }$.

Note that $Y=G_1/K_1$ is a symmetric space of nonpositive curvature,
with $G_1\subseteq G$ and $K_1\subseteq K$. Since $Y$ is the fixed
point set of $d$--isometries, so it is $d$--convex, and $G_1$ acting
transitively on $Y$ as both $d^{\ast }$ and $d$--isometries.

By Sublemma, either $Y=F$, in which case $F$ is a $d$--flat, or
$Y=({\mathbb H}^2)^r$. In this case, the arguments in the previous
subsections implies that $(Y,d)$ must split as the product of $r$
surfaces, each is symmetric since $({\mathbb H}^2, d^{\ast })$ has
rank $1$. So $d$ is a symmetric metric on $Y$, and $F$ is a $d$--flat.

Now consider any maximal flat $F\subseteq B$ passing through
$x_0$. The restriction of $d$ to each flat $F$ is a Euclidean metric
which is invariant under the Weyl group, which acts irreducibly on $F$
since $G$ is assumed to be simple. So $d|_F=\lambda (F)d^{\ast }|_F$.
Since any two flats through $x_0$ can be joined by finitely many
flats through $x_0$ so that each one intersects the next at more than
one point,  $\lambda (F)$ must be constant, and $d=\lambda d^{\ast
}$ is symmetric. This completes the proof of Lemma~\ref{G-inv}
assuming Sublemma. \end{proof}

\begin{proof}[Proof of Sublemma] Let us start with part (1). Without
loss of generality, we may assume that $G$ is simple. Assume $M=1$, we
want to conclude that $G/K$ is the hyperbolic plane. Let $r$ be the
rank of $G/K$.

Since the Lie algebra $ \mathfrak{ m}$ of $M$ is trivial, $\mathfrak{
  h}_{\mathbb C}$ is a Cartan subalgebra for $\mathfrak{ g}_{\mathbb
  C}$, the complexification of $\mathfrak{g}$. Denote by $G'$ the real
  points of $G_{\mathbb C}$ under the conjugation of $\mathfrak{g}$ in
  $\mathfrak{ g}_{\mathbb C}$. Then $G = \mbox{Ad}(G)$ is the identity
  component of $G'$.  Denote by $K'$ the maximal compact subgroup of
  $G'$ with identity component $K$.

Let $P\cong ({\mathbb Z}_2)^r$ be the group of $2$--torsion elements in
$H_{\mathbb C}=\mbox{exp} (\mathfrak{ h}_{\mathbb C})$. Then $P$ is
contained in $Z_{K'}(A)$, the centralizer of $A$ in $K'$.

Note that the number of connected components in $G'$ is bounded by the
order of the outer automorphism group $Q=Aut(G)/Int(G)$. So the
assumption $M=Z_{K}(A)=1$ implies that $2^r \leq |Q|$.

By the Table 10 on page 156 of \cite{L}, we see that this is
impossible unless $\mathfrak{g}_{\mathbb C}$ is $\mathfrak{ a}_1=
\mathfrak{ sl}(2,{\mathbb C})$ or $\mathfrak{ d}_4= \mathfrak{
so}(8,{\mathbb C})$. In the $\mathfrak{ d}_4$ case, in fact in all
four classical Lie algebra cases, it can be easily checked that $M\neq
1$ except for $\mathfrak{a}_1$.  This completes the proof of part (1).

For part (2), let us apply part (1) to the symmetric space (which may
contain a Euclidean factor now) $Y$, since $M_1=Z_{K_1}(A)$ is
trivial, by part (1), we know that $Y=({\mathbb H}^2)^s \times
{\mathbb R}^{r-s}$. We want to show that $s$ is either $0$ or $r$ if
$G$ is simple.

Let $M'=N_K(A)$ be the normalizer of $A$ in $K$. Then $M'$ ( or $M$)
consists of elements in $K$ that stabilize $F$ ( fix every point in
$F$). $M$ is a normal subgroup of finite index in $M'$, and $W=M'/M$
is Weyl group which acts irreducibly on $F$. For each $g\in M'$, since
$gMg^{-1}=M$, we have $g(Y)=Y$.  Now since any isometry of $Y$
preserves the Euclidean factor, we know that the ${\mathbb R}^{r-s}$
part of $F$ is invariant under $W$. So either $s=0$ or $s=r$. This
completes the proof of Sublemma.
\end{proof}

\begin{Remark} Note that the CAT(0) assumption on $d$ in
Lemma~\ref{G-inv} is necessary. When rank $r\geq 2$, there are
$G$--invariant length metric on $G/K$ which are  geodesic. The
simplest example would be taking product of two copies of symmetric
spaces, and let $d=d_1+d_2$ where $d_i$ are the standard metrics on
the factors. Then clearly $d$ is $G$--invariant but not Riemannian.
More generally, Planche showed \cite{P} that the set of
$G$--invariant Finsler metrics on $G/K$ is in one-to-one
correspondence with the set of $W$--invariant norms on a maximal flat
$F$. (Here $W$ denotes the Weyl group.)  Some of these metrics are even 
{\em uniquely\/} geodesic, that is, any two points in $G/K$ can be joined by a unique geodesic segment.
\end{Remark}\np


\begin{thebibliography}

\bibitem{B} {\bf W Ballmann}, {\it Lectures on Spaces of Nonpositive
    Curvature,} Birkh\"{a}user (1995)

\bibitem{BGS} {\bf W Ballmann}, {\bf M Gromov}, {\bf V
Schroeder},
 {\it Manifolds of Nonpositive Curvature}, Birkh\"{a}user (1985)

\bibitem{BH} {\bf M Bridson}, {\bf A Haefliger}, {\it Spaces of
nonpositive curvature}, Springer--Verlag, Berlin, Heidelberg, New York
(1999)


\bibitem{CE} {\bf S-S Chen}, {\bf P Eberlein}, {\it Isometry groups of
simply connected manifolds of non-positive curvature}, Illinois
Journal of Math.  24 (1980) 73--103

\bibitem{CD1} {\bf R Charney}, {\bf M\,W Davis}, {\it Singular metrics
    of nonpositive curvature on branched covers of Riemannian
    manifolds}, Amer.  J. Math.  115 (1993) 929--1009

\bibitem{CD2} {\bf R Charney}, {\bf M\,W Davis}, {\it The polar dual of
    a convex polyhedral set in hyperbolic space}, Michigan Math. J.
     42 (1995) 479--510

\bibitem{CDM} {\bf R Charney}, {\bf M\,W Davis}, {\bf G Moussong}, {\it
    Nonpositively curved piecewise Euclidean structures on hyperbolic
    manifolds}, Michigan Math. J.   44 (1997) 201--208

\bibitem{DJ} {\bf M\,W Davis}, {\bf T Januszkiewicz}, {\it
Hyperbolization of polyhedra}, J. Diff. Geom. 34 (1991)
347--388


\bibitem{E1} {\bf P Eberlein}, {\it Rigidity of lattices of
    nonpositive curvature}, Ergodic Theory \& Dynamical Systems,
    3 (1983) 47--85

\bibitem{E2} {\bf P Eberlein}, {\it Isometry groups of simply connected
    manifolds of non-positive curvature, II}, Acta Math.  149
    (1982) 41--69

\bibitem{E3} {\bf P Eberlein}, {\it Euclidean de Rham factor of a
lattice of nonpositive curvature}, J. Diff. Geometry, 18 (1983) 209--220

\bibitem{E4} {\bf P Eberlein}, {\it Lattices in spaces of nonpositive
    curvature}, Ann. Math. 111 (1980) 435--476

\bibitem{EF} {\bf A Eskin}, {\bf B Farb}, {\it Quasi-flats and rigidity
    in higher rank symmetric spaces}, J. Amer. Math. Soc. 10 (1997) 653--692

\bibitem{G} {\bf M Gromov}, {\it Hyperbolic groups}, Essays in group
theory (S\,M Gersten, ed.), MSRI Publ. 8, Springer--Verlag,
Berlin, Heidelberg, New York (1987)  75--264

\bibitem{GT} {\bf M Gromov}, {\bf W Thurston}, {\it Pinching constants
for hyperbolic manifolds}, Invent. Math.  89 (1987)  1--12

\bibitem{KL} {\bf B Kleiner}, {\bf B Leeb}, {\it Rigidity of
quasi-isometries for symmetric spaces and Euclidean buildings},
 Inst. Hautes \'{E}tudes Sci. Publ. Math. No. 86 (1997) 115--197 
  
\bibitem{Le}{\bf  B Leeb}, {\it A characterization of irreducible
symmetric spaces and Euclidean buildings of higher rank by their
asymptotic geometry}, preprint (1997)

\bibitem{L}{\bf  O Loos}, {\it Symmetric spaces: Volume II}, W.A.
Benjamin, Inc. (1969)

\bibitem{M} {\bf G Mostow}, {\it Strong rigidity of locally symmetric
spaces}, Annals of Math. Studies 78, Princeton University Press (1973)

\bibitem{MZ} {\bf D Montgomery}, {\bf  L Zippin}, {\it Topological
    Transformation Groups}, Interscience, New York (1955)

\bibitem{P}{\bf  P Planche}, {\it Structures de Finsler invariantes sur
les espaces sym\'{e}triques}, CR. Acad. Sci. Paris,  321 (1995) 1455--1458

\bibitem{R}{\bf  M\,S Raghunathan}, {\it Discrete subgroups of Lie
groups}, Springer--Verlag, Berlin, Heidelberg, New York (1972)


\bibitem{S}{\bf  A Selberg}, {\it Recent developments in the theory of
discontinuous groups of motion of symmetric spaces}, Springer Lecture
Notes in Math. 118, Springer--Verlag, Berlin, Heidelberg, New York (1969)

\end{thebibliography}
\end{document}